\documentclass[10pt]{article}
\font\smallit=cmti10

\usepackage{tabularx}
\usepackage{amssymb,amsmath,amsthm,enumitem} 
\usepackage{subcaption,booktabs}
\usepackage{hyperref}
\makeatletter

\renewcommand\section{\@startsection {section}{1}{\z@}
{-30pt \@plus -1ex \@minus -.2ex}
{2.3ex \@plus.2ex}
{\normalfont\normalsize\bfseries\boldmath}}

\renewcommand\subsection{\@startsection{subsection}{2}{\z@}
{-3.25ex\@plus -1ex \@minus -.2ex}
{1.5ex \@plus .2ex}
{\normalfont\normalsize\bfseries\boldmath}}

\renewcommand{\@seccntformat}[1]{\csname the#1\endcsname. }

\makeatother

\newtheorem{theorem}{Theorem}
\newtheorem{lemma}{Lemma}

\newtheorem{proposition}{Proposition}

\newtheorem{conjecture}{Conjecture}
\newtheorem{observation}{Observation}
\theoremstyle{definition}

\newtheorem{remark}{Remark}

\def\R{\mathcal{R}}
\def\M{\mathcal{M}}

\begin{document}
\begin{center}
\uppercase{Screw discrete dynamical systems and their 
applications to exact slow NIM} 

\vskip 20pt
{\bf Vladimir Gurvich}\\
{\smallit National Research University Higher School of Economics (HSE), Moscow, Russia}\\
{\tt vgurvich@hse.ru}, {\tt vladimir.gurvich@gmail.com}\\
\vskip 10pt
{\bf Mariya Naumova}\\
{\smallit Rutgers Business School, Rutgers University, Piscataway, NJ, United States}\\
{\tt mnaumova@business.rutgers.edu}\\
\end{center}
\vskip 20pt
\centerline{\smallit Received: , Revised: , Accepted: , Published: } 
\vskip 30pt

\centerline{\bf Abstract}
\noindent
Given integers  $n,k,\ell$  such that  $0<k<n, \; 1<\ell$    
and an integer vector 
\newline 
$x = (x_1,\ldots,x_n)$; 
denote by $m = m(x)$  the number of entries of  $x$  
that are multiple of  $\ell$. 
\newline 
Choose  $n-k$  entries of $x$ as follows:
if  $n-k \leq m(x)$,  take  $n-k$  
smallest entries of $x$  multiple of $\ell$; 
if  $n-k > m(x)$, take all $m$ such entries,  
if any, and add remaining  $n-k-m$  entries arbitrarily, 
for example, take the largest ones. 
\newline
In one step, the chosen  $n-k$  entries 
(bears) keep their values, 
while the remaining  $k$  (bulls) are reduced by 1. 
Repeat such steps getting the sequence  
\newline 
$S = S(n,k,\ell,x^0) = (x^0 \to x^1 \to \ldots \to x^j \to \ldots)$.  
It is ``quasi-periodic". 
\newline 
More precisely, 
there is a function $N = N(n,k,\ell,x^0)$  such that 
for all $j \geq N$ we have 
$m(x^j) \geq n-k$ and $range(x^j) \leq \ell$, where 
$range(x) = (\max(x_i \mid i \in [n]) - \min(x_i \mid i \in [n])$.
Furthermore, $N$ is a polynomial in $n,k,\ell,$ and $range(x^0)$ 
and can be computed in time linear in 
$n,k,\ell$, and $\log(1 + range(x^0))$.
\newline 
After $N$ steps, the system  moves ``like a screw". 
Assuming that $x_1 \leq \dots \leq x_n$, 
introduce the cyclical order on $[n] = \{1, \ldots, n\}$ 
considering  1 and $n$ as neighbors. 
Then, bears and bulls partition $[n]$ into two intervals, 
rotating by the angle $2 \pi k /n$  with every $\ell$  steps. 
Furthermore, after every 
$p = \ell n / GCD(n,k) = \ell LCM(n,k) / k$  steps  
all entries of $x$  are reduced by the same value  $\delta = pk/n$, 
that is, $x_i^{j+p} - x_i^j = \delta$ 
for all $i \in [n]$ and  $j \geq N$. 
We provide an algorithm computing  $N$ 
(and  $x^j$)  in time linear in  
$n,k,\ell, \log(1 + range(x^0))$ (and $\log (1+j)$).  
\newline 
In case  $k=n-1$ and $\ell = 2$  
such screw dynamical system are applicable to impartial games; 
see  Gurvich, Martynov, Maximchuk, and Vyalyi,  
``On Remoteness Functions of Exact Slow $k$-NIM with $k+1$ Piles",  
https://arxiv.org/abs/2304.06498 (2023).
\newline 
\noindent{\bf AMS subjects:  91A05, 91A46, 91A68} 

\section{Introduction} 
\label{s0}
\subsection{Main concepts}
\label{ss00}
Fix integers $n,k,\ell$ such that $0 < k < n, \; 1 < \ell$, 
and consider an integer vector $x = (x_1, \ldots, x_n)$  
(that may have negative entries).   
We assume that $x$ is defined up to a permutation of its entries. 
Hence, without loss of generality (wlog),  
we will require monotonicity: $x_1 \leq \ldots \leq x_n$. 

Denote by $m$ the number of entries of $x$  
that are multiple of  $\ell$. 
If  $m \geq n-k$, choose the smallest $n-k$ of such entries.
If  $m < n-k$, choose 
all  $m$  such entries 
and add arbitrary $n-k-m$  others, 
for example the largest ones. 
Note that in both cases the choice may be not unique, 
since entries of $x$  may be equal.
In case of such ambiguity we use the following 

\medskip 
\noindent 
{\em tie-breaking} rule: 
always choose the rightmost entries, 
that is, $x_i$  with the largest $i$. 

\medskip 

The chosen $n-k$ entries of  $x$  
will be called {\em bearish} and 
their indices - {\em bears}.  
 
For any feasible $n,k,\ell$, we define a discrete dynamical system 
by the following {\em GM-$(n,k,\ell)$-rule}. 
Given an integer vector $x$, 
we define a move $x \to x'$:   
the $n-k$ bearish entries of $x$ do not change, 
while the remaining $k$, bullish, entries are reduced by 1.
This (unique) move $x \to x'$  and 
the (unique) sequence 
$$S = S(n,k,\ell,x) = (x = x^0 \to x^1 \to \ldots \to x^j \to \ldots)$$   
will be called  {\em GM-$(n,k,\ell)$-move} and 
{\em GM-$(n,k,\ell)$-sequence}, respectively. 
We will omit arguments 
$(n,k,\ell)$ whenever their values are unambiguous from the context 
and simply say GM-rule, GM-move, and GM-sequence, for short. 

\subsection{Preliminary results}
\label{ss01}
 
We begin with several simple properties 
of GM-moves and GM-sequences. 

\begin{observation}
\label{o1}
A GM-move $x \to x'$ 
respects monotonicity, that is,   
$x'_1 \leq \ldots \leq x'_n$ whenever 
$x_1 \leq \ldots \leq x_n$.  
\end{observation}

\proof 
It follows immediately from the tie-breaking rule. 
\qed 

\medskip 

Recall that $m = m(x)$ denotes the number of entries in $x$  
that are multiple of  $\ell$. 

\begin{observation}
\label{o2}
A GM-$(n,k,\ell)$-move $x \to x'$ 
respects inequality  $n-k \leq m(x)$,  
that is, $n-k \leq m(x')$  whenever $n-k \leq m(x)$.  
\end{observation}

\proof 
Assume that $i \in [n] = \{1,\ldots,n\}$ 
``de-bears"  a move $x \to x'$, 
that is, $i$ is a bear in $x$  but is a bull in $x'$. 
Then, $x_i$  is a multiple of $\ell$ and $x'_i = x_i$, 
since $i$ is a bear in $x$. 
In contrast, $i$  is a bull in $x'$.
Hence, by the tie-breaking rule, 
there exist $n-k$ bears of $x'$ 
each of which is either smaller than $i$, 
or larger than $i$, but its 
bearish value equals  $x_i = x'_i$. 
Thus, $i$  may stop to be a bear,  
after a GM-move $x \to x'$,  
only if a ``replacing bear" appears in  $x'$.
\qed 

\begin{observation}
\label{o3}
Given an arbitrary GM-sequence $S(n,k,\ell,x^0)$,
if inequality  $n-k \leq m(x^j)$ fails for $j=0$ 
then it holds for all  $j \geq \frac{1}{k} 
[\sum_{i=1}^{m(x^0)-(n-k)} (x_i^0 \mod \ell)] - m(x^0)$. 
\end{observation} 

In other words, inequality $n-k \leq m(x^j)$ 
holds for all GM-sequences $S(n,k,\ell,x^0)$  
and  $j \geq \frac{1}{k} 
[\sum_{i=1}^{m(x^0)-(n-k)} (x_i^0 \mod \ell)] - m(x^0)$, 
where, by convention, the sum equals 0 when $m(x^0)-(n-k) \leq 0$.

\proof 
By Observation \ref{o2}, 
this claim holds for all $j \geq 0$ if  $n-k \leq m(x^0)$.
Suppose that  $n-k > m(x^j)$  for $j = 0, \dots, j_0$.
Then,  for such  $j$  we have 
$m(x^j) \geq 0$  entries $x_i^j$  
that are bearish,  multiple of $\ell$, and 
are not changed by a GM-move $x^j \to x^{j+1}$. 

Also  $x^j$  has exactly $k$  bullish entries, 
which are not multiple of $\ell$  
and reduced by 1 by each GM-move,  
until inequality  $n-k \leq m(x^j)$ fails.
Hence, it will be achieved 
(and then hold, by Observation \ref{o2}) 
in at most 
$\frac{1}{k} [\sum_{i=1}^{n-k-m(x^0)} (x_i^0 \mod \ell)] - m(x^0)$ 
GM-moves from $x^0$. 
Either this number is positive, 
or $n-k \leq m(x^j)$ holds already for $j=0$. 
\qed 

\medskip 

Given $x = (x_i \mid i \in [n]) = \{1, \ldots, n\}$, let us set 
$$range(x) = \max(x_i \mid i \in [n]) - \min(x_i \mid i \in [n]) 
= x_n - x_1.$$ 

\begin{observation}
\label{o4}
One move $x \to x'$ may either keep the range, $range(x') = range(x)$,  
or change it by at most 1, $range(x') = range(x) \pm 1$.   
A GM-move $x \to x'$ respects the inequality 
$range(x) \leq \ell$. 
In other words, 
$range(x') \leq \ell$ whenever $range(x) \leq \ell$. 
\end{observation} 

\proof 
The first claim is obvious. 
The second one is not difficult too. 
Assume for contradiction that 
$range(x) = \ell$, while $range(x') = \ell+1$. 
for a GM-move $x \to x'$. 
This can happen only if 
$x'_1 = x_1 - 1$, while  $x'_n = x_n$. 
The second equality holds 
if and only if  $n$  is a bear in $x$. 
Hence, $x_n = a \ell$  for some integer  $a$. 
Furthermore, $range(x) = \ell$ and, hence, $x_1 = (a-1) \ell$. 
Then, by the tie-breaking rule, 
$1$  is a bear in  $x$  and $x_1 = x'_1$ 
unless $x_1 = \ldots = x_{n-k} = x_{n-k+1} = (a-1) \ell$.  
Yet, in this case $n$ is not a bear in $x$ 
and $x'_n = a \ell-1$, which is a contradiction. 
\qed 

\medskip 

Thus, for each GM-sequence  $S = S(n,k,\ell,x^0)$  
there exists a unique integer   
$N = N(n,k,\ell,x^0) \geq 0$  such that  
$m(x^j) \geq n-k$  and  
$range(x^j) = x_n^j - x_1^j \leq \ell$ 
if and only if  $j \geq N$. 
In other words, $S$  
is partitioned into two intervals: an initial one  
$S_0 = (x_0 \to x^1 \to \ldots \to x^{N-1})$ 
and an infinite one   
$S_\infty = (x^N \to x^{N+1} \to \dots)$. 
By convention, we set $N=0$ and $S_0 = \emptyset$ if 
both inequalities hold already for $x^0$. 

The value of $m(x^j) - (n-k)$ is negative 
and monotone non-decreasing for $j < N$ and 
it becomes and remains nonnegative for $j \geq N$. 
Note, however, that it may decrease 
with a move $x^j \to x^{j+1}$  when $j \geq N$, 
because several bears in $x^j$
may disappear after this move if  $m(x^j) \geq n-k$; 
see examples of Section \ref{s-examples}. 


\subsection{Main results}
\label{ss02} 
We proceed with more complicated statements, 
which will be proven later. 

\begin{proposition}
\label{p1}
For a GM-sequence $S(n,k,\ell,x) = S = S_0 \cup S_\infty$, 
the inequality $range(x^j) \leq \ell$ holds 
if and only if  $j \in S_\infty$, that is, $j \geq N$.
Furthermore $N = N(n,k,\ell,x)$  
is a polynomial in its arguments and 
can be computed in time linear in $n,k,\ell$, 
and $\log(1 + range(x^0)) = \log(1 + x^0_n - x^0_1)$. 
\end{proposition}

Define a period $p$ and a reduction value  $\delta$  by formulas:
\begin{equation} 
\label{eq-main}
p = p(n,k,\ell) = \ell \, n / GCD(n,k) = \ell \; LCM(n,k)/k, 
\;\; \delta = \delta(n,k,\ell) = pk/n.
\end{equation} 
Here  GCD and LCM stand for 
the {\em greatest common divisor} and {\em least common multiple}. 
It is both obvious and well known that  $GCD(n,k) LCM(n,k) = nk$.
Furthermore, $GCD(n,k) = 1$ if and only if  $LCM(n,k) = nk$.
In this case $n$ and $k$ are called {\em co-prime}  and  
we have  $p = \ell n$ and  $\delta = \ell k$.  
For example, this is the case if  $k = 1$ or $k = n-1$. 
Our main result is given by the following statement. 

\begin{theorem}
\label{t-main} 
For every $i \in [n]$ and $j \geq N$ we have  
 $x_i^{j+p} - x_i^j = \delta$.  
\end{theorem}

Note that  $p$  is the ``minimal period" 
if  $k=1$  or $k=n-1$. 
If  $1 < k < n-1$ then equation  
$x_i^{j+q} - x_i^j = \delta$  may hold      
for every $i \in [n]$ and $j \geq N$  
already for a divisor  $q$  of  $p$. 
Such an example will be given in Section \ref{s-examples}. 

Note also that if $x_i^{j+p} - x_i^j = \delta_j$ 
and $\delta_1 = \ldots = \delta_n = \delta$ 
then $\delta = pk/n$.
Indeed, each of $p$ GM-moves reduces exactly $k$ 
of $n$ entries and each of them - exactly by 1. 

\begin{remark}
The periodical phase $S_\infty$ 
of the considered discrete dynamical system   
is similar to a screw. 
Introduce a cyclical structure on 
$[n] = \{1, \ldots, n\}$ treating 1 and $n$ as neighbors. 
Cycle  $[n]$  is partitioned into 
two intervals, of size $k$ (bulls) and $n-k$ (bears). 
Every $\ell$ GM-moves rotate 
these two intervals by the angle  $2 \pi k / n$.  
Furthermore, with every  $p$  GM-moves 
this screw advances forward by $\delta$.  
\end{remark}

\medskip 

The length of a maximal interval of successive bulls in $S_\infty$ 
is a multiple of $\ell$. 
More precisely, the following long but simple statement holds. 

\begin{observation}
\label{04} 
Fix an $i \in [n]$ and consider the interval     
$$J = J(i) = \{j_0-1, j_0, j_0+1, \ldots, j_0+t-1, j_0+t\} \subset S_\infty$$
such that $i$ is a bear in two vectors $x^j$,   
for $j = j_0-1, \; j=j_0+t$,  
and a bull the remaining $t$ vectors, with 
$j \in \{j_0, \ldots, j+t-1\}$. 
Then  $t$  is a multiple of $\ell$. 

Also  $x_i^{j_0-1} = x_i^{j_0}$  and 
$x_i^{j_0+t}$ are multiples of $\ell$, 
while $x_i^{j_0}, \ldots, x_i^{j_0+t-1}$ are not.  
\end{observation}

\proof 
If  $j \in S_\infty$  then 
$range(x^j) \leq \ell$  and  $m(x^j) \geq n-k$. 
Thus, $x_i^j$ is a multiple of $\ell$ 
whenever $i$ is a bear in $x^j$, 
but not vice versa. 
In particular, $x_i^{j_0-1}$ 
is a multiple of $\ell$, since $x_i^{j_0-1}$ 
is a bearish entry,  while 
$x_i^{j_0}$  is also a multiple of $\ell$, 
since  $x_i^{j_0} = x_i^{j_0}$, but it is bullish.  
Furthermore, $x_i^{j_0+t}$ is bearish and,   
hence, a multiple of $\ell$ again.  
According to the GM-rule, 
the entries between the considered two 
are reduced by 1 with each GM-move. 
Thus, $t$ is a multiple of $\ell$ too. 
\qed

\medskip 

We assumed that $j_0-1 \in S_\infty$. 
Yet, this assumption can be waved 
if we consider  $j \mod p$  rather than $j$.  
The length of a maximal interval 
of successive bulls in $S_\infty$ 
is a multiple of $\ell$  
if we assume that  $1$ and $p$ are neighbors, 
$p+1 = 1 \mod p$. 


\subsection{Plan of the proof of Theorem \ref{t-main}} 
Introduce the cyclical order over $[n] = \{1, \ldots, n\}$. 
By convention, 1 and $n$ are neighbors.  
We will show that some properties of GM-sequences can be 
formulated in terms of this cyclical order, 
which looks somewhat surprising in presence of monotonicity 
$x_1 \leq \ldots \leq x_n$.  

A set of successive elements of $[n]$ 
in the cyclical order is called an {\em interval}. 
By convention, a single-element set  
and even the empty set are intervals too. 
We say that an interval is {\em of type 2}  
if it contains both  1 and $n$, 
otherwise it is {\em of type 1}.

\begin{observation}
Partition $[n]$ into two complementary subsets  
$I$  and $I^c = [n] \setminus I$. 
\begin{itemize}
\item[(j)] 
$I$ is an interval 
if and only if $I^c$ is an interval.
\item[(jj)] At most one of these two intervals is of type 2.
\item[(jjj)] Both are of type 1 if and only if 
one contains 1, while the other contains  $n$.
\item[(jv)] Any subinterval of an interval of type 1 
is of type 1, too. 
\end{itemize}
\end{observation}

\proof It is straightforward. 
\qed 

\medskip 

Given, $n,k,\ell$, and $x$, 
denote by $M = M(x)$ the set of entries of 
$x$ that are multiples of $\ell$ and 
by $P = P(x)$ the set of bears in $x$. 
Recall that $|M| = m$.  
By the GM-rule, $|P| = n-k$  and     
$P \subseteq M$  if and only if  $n-k \leq m$.

\begin{proposition}
\label{p3}
If $n-k \leq m$ and $range(x) \leq \ell$ 
then both sets $M$ and $P$ are intervals of $[n]$. 
\end{proposition}

\proof 
By our notation and assumptions, 
$n \geq |M| = m \geq n-k$. 
Since $range(x) \leq \ell$, 
the $m$ entries of $x$ that are multiple of $\ell$    
may take either (i) the same value $a \ell$,    
or (ii) two values $a \ell$  and  $(a+1) \ell$, 
for some integer  $a$. 

Case (i). If  $m=n$  then $M = [n]$   
is an interval of type 2, 
while  $P = \{k+1, \ldots, n\}$  
is an interval of type 1, since $k > 0$.  
If $m < n$  then $P \subseteq M \subset [n]$ 
and, by the tie-breaking rule, $P$  is the rightmost 
subinterval of $M$  of size $n-k$. 
In both subcases $P$  is of type 1,  
since it cannot contain both 1 and $n$, 
because $|P| = n-k < n$. 

\medskip 

Case (ii). Clearly,   
$x_1 = a \ell, \; x_n = (a+1) \ell$, and  
set $M$ is partitioned into two non-empty subintervals:  
$M = M_1 \cup M_n$, where $1 \in M_1$  and  $n \in M_n$. 
Thus, $M$ is of type 2, since it contains both $1$ and $n$. 

By the GM-rule, $P \cap M_1 \neq \emptyset$. 
There are two subcases: 
either  $|M_1| \geq n-k$, 
then  $P \subseteq M_1$, or 
$|M_1| < n-k$, then  
($P \supseteq M_1$  and  $P \cap M_n \neq \emptyset$). 

In the first subcase, by the  tie-breaking rule,  
$P$  is the rightmost subinterval of  $M_1$  of size $n-k$ 
and, hence, $P$ is of type 1 (while $M$ is of type 2).

In the second subcase, $P = M_1 \cup M'_n$, 
where $M'_n$ is the rightmost subinterval of $M_n$ 
of size  $n-k-|M_1|$. 
In this case $P$  (as well as $M$) is 
an interval of type 2. 
\qed 

\medskip 

We will derive Theorem \ref{t-main}, 
from the following main lemma. 

\begin{lemma}
\label{l-main} 
Fix feasible  $n,k,\ell,$  
$x^0 = (x^0_1, \ldots, x^0_n)$,    
and consider the (unique) 
GM-sequence  $S = S(n,k,\ell,x^0)$.   
Then, $i \in [n]$  is a bear in  $x^j$  
if and only if 
$i + k \mod n$  is a bear in  $x^{j + \ell}$,  
for each $i \in [n]$ and $j \geq N(n,k,\ell,x^0)$.   
\end{lemma}

In other words, for each  $x \in S_\infty$, 
in every $\ell$ GM-moves from  $x$, both intervals, 
of $n-k$  bears and  $k$  bulls,   
are shifted to the right by  $k$, 
or equivalently, to the left, by $n-k$, 
in the cyclical order over $[n]$. 

For example, if  $k = n-1$,  
the unique bear is shifted by 1 to the left, 
modulo $n$. 
In particular, 1 is shifted to $n$, 
since $1-1 = 0 = n \mod n$.

If  $k = 1$, then the interval 
$[n] \setminus \{i\}$  of $n-1$ bears 
is shifted by 1 to the right, modulo $n$. 
In particular, $[n] \setminus \{n\}$ 
is shifted to  $[n] \setminus \{1\}$,  
since $n+1 = 1 \mod n$. 

In general, in every $p$  GM-moves 
each  $i \in [n]$  will be a bear 
the same number of times, $\delta$, 
where $p$ and $\delta$ are defined by 
formula (\ref{eq-main}). 
Hence, in every $p$ GM-moves from $x^j \in S_\infty$,  
each entry $x^j_i$  will be reduced by $\delta$, 
that is, $x_i^{j+p} - x_i^j = \delta$  for every 
$i \in [n]$ and $j \geq N(n,k,\ell,x^0)$, 
implying Theorem~\ref{t-main}.

\subsection{Applications to impartial games} 
Sprague \cite{Spr36} and Grundy \cite{Gru39} 
introduced the Sprague-Grundy function 
of impartial games. 
Then, Smith \cite{Smi66} defined  
for them the remoteness function. 
We will apply the GM-$(n,k,\ell)$-sequences 
with $\ell = 2$  for computing 
the remoteness function 
of the game NIM$(n,k)$ of Exact Slow NIM. 
In this game two players alternate turns 
reducing by each move any $k$ from $n$ piles of stones,   
by one stone each pile. 
See the last section for definitions and more details. 

\section{Examples illustrating main lemma}
\label{s-examples}
We will represent a GM-sequence 
$x^0 \to x^1 \to \ldots \to x^j \to \ldots$ 
by a table whose rows and columns 
are numbered by $j = 0,1, \ldots$ 
and  $i = 1, \ldots, n,$ respectively. 
The following four examples illustrate Lemma \ref{l-main} 
for $\ell = k = n-1 = 3$. 
In every row there exists a unique bear, since $n-k = 4-3 =1$. 
In each 3 GM-moves, this bear is shifted to the left 
(that is, reduced) by 1 in the cyclical order. 
In particular, 1 is shifted to $n=4$. 

\begin{table}[h] 
\setlength\tabcolsep{0pt} 
\small 
\caption{Examples to Lemma \ref{l-main} 
for $n=4, k = 3$ and  $\ell = 3$.} 
\label{tab:4subt}
\begin{subtable}{0.22\textwidth}
\begin{tabular*}{\linewidth}{@{\extracolsep{\fill}}ccccc}
\toprule
$x_1$ & $x_2$ & $x_3$ & $x_4$  &  bears \\
\midrule
15& 15& 17& 18 & 2  \\
14& 15& 16& 17 & 2  \\
13& 15& 15& 16 & 3  \\
\hline
12& 14& 15& 15 & 1  \\
12& 13& 14& 14 & 1  \\
12& 12& 13& 13 & 2  \\
\hline
11 &12& 12& 12 & 4  \\
10& 11& 11& 12 & 4 \\
9& 10 &10 &12  & 1 \\
\hline
9& 9& 9& 11 & 3  \\
8& 8& 9& 10 & 3 \\
7& 7& 9& 9  & 4 \\
\hline
6& 6& 8 &9 &  2  \\
\bottomrule
\end{tabular*}
\label{tab:sub1}
\end{subtable}
\hfill
\begin{subtable}{0.22\textwidth}
\begin{tabular*}{\linewidth}{@{\extracolsep{\fill}}ccccc}
\toprule
$x_1$ & $x_2$ & $x_3$ & $x_4$  &  bears\\
\midrule
15& 16& 17& 17 & 1 \\
15& 15& 16& 16 & 2 \\
14& 15& 15& 15 & 4 \\ 
\hline
13& 14& 14& 15& 4 \\
12& 13& 13& 15& 1\\ 
12& 12& 12& 14& 3\\
\hline
11& 11& 12& 13& 3\\  
10& 10& 12& 13& 4\\
9 & 9& 11& 12& 2\\
\hline
8 &9 &10 &11 & 2\\ 
7 &9 &9 &10 & 3\\
6 &8 &9 &9  &1\\
\hline
6 &7 &8 &8 &1 \\  
\bottomrule
\end{tabular*}
\label{tab:sub2}
\end{subtable}%
\hfill
\begin{subtable}{0.22\textwidth}
\begin{tabular*}{\linewidth}{@{\extracolsep{\fill}}ccccc}
\toprule
$x_1$ & $x_2$ & $x_3$ & $x_4$  &  bears \\
\midrule
15& 17& 17& 18&1 \\
15& 16 &16 &17&1 \\
15& 15& 15& 16 & 3 \\
\hline
14& 14& 15 &15& 4 \\
13& 13& 14 &15&  4 \\
12& 12& 13& 15 & 2 \\
\hline
11 &12& 12& 14&3  \\
10& 11& 12& 13 &3\\
9 & 10& 12& 12& 1\\
\hline
9 &9& 11 &11& 2 \\
9 &9 &10& 10 &2\\
7 &9& 9& 9 & 4\\
\hline
6& 8 &8& 9& 1\\
\bottomrule
\end{tabular*}
 \label{tab:sub3}
\end{subtable}
\hfill
\begin{subtable}{0.22\textwidth}
\begin{tabular*}{\linewidth}{@{\extracolsep{\fill}}ccccc} 
\toprule
$x_1$ & $x_2$ & $x_3$ & $x_4$  &  bears \\
\midrule
15 &  18 &  18  & 18  &  1  \\
15 & 17  & 17 & 17 &  1  \\
15 &  16 &  16  & 16  &  1  \\
\hline
15 &  15  & 15 &  15  &   4  \\
14  & 14  & 14  & 15 &   4  \\
13  & 13  & 13 & 15  &   4  \\
\hline
12  & 12 &  12 &  15 &  3   \\
11  & 11 &  12  & 14 &  3 \\
10 &  10  & 12 &  13   &  3 \\
\hline
9  & 9 &  12 &  12  &   2 \\ 
8  & 9 &  11 &  11  &   2 \\
7  & 9  & 10 &  10  &  2 \\
\hline
6  & 9  & 9 &  9 & 1  \\ 
\bottomrule
\end{tabular*}
\label{tab:sub4}
\end{subtable}%
\hfill
\end{table}

Table 2 shows 4 more examples with $n=5, \ell = 3, \; k = 1,2,3,4$. 

Since $n = 5$ is prime, $n$ and $k$ are co-prime, hence, $p = n \ell = 15$ and $\delta = pk/n = 9$ for all $k$.   
Yet, for $k=3$ there exists a smaller period 
$q = 5$  with $\delta = q k / n = 3$.

\begin{table}[t] 
\setlength\tabcolsep{0pt}
\small 
\caption{Examples to Lemma \ref{l-main} 
for $n=5, k = 1, 2, 3, 4$, and  $\ell = 3$.} 
\begin{subtable}{0.28\textwidth}
\begin{tabular*}{\linewidth}[t]{@{\extracolsep{\fill}}ccccc|c}
\toprule
$x_1$ & $x_2$ & $x_3$ & $x_4$  & $x_5$ & bears\\
\midrule
3& 3& 4& 6& 6& $[5]\setminus\{3\}$\\
3& 3& 3& 6& 6& $[5] \setminus \{4\}$\\
3& 3& 3& 5& 6& $[5] \setminus \{4\}$\\
\hline
3& 3& 3& 4& 6& $[5] \setminus \{4\}$\\
3& 3 &3& 3& 6 & $[5] \setminus \{5\}$\\
3& 3& 3& 3& 5& $[5] \setminus \{5\}$\\
\hline
3& 3& 3& 3& 4& $[5] \setminus \{5\}$\\
3& 3& 3& 3& 3& $[5] \setminus \{1\}$\\
2& 3& 3& 3& 3& $[5] \setminus \{1\}$\\
\hline
1& 3& 3& 3& 3& $[5] \setminus \{1\}$\\
0& 3& 3& 3& 3& $[5] \setminus \{2\}$\\
0& 2& 3& 3& 3& $[5] \setminus \{2\}$\\
\hline
0& 1& 3& 3& 3& $[5] \setminus \{2\}$\\
0& 0& 3& 3& 3& $[5] \setminus \{3\}$\\
0& 0& 2& 3& 3& $[5] \setminus \{3\}$\\
\hline
0& 0& 1& 3& 3& $[5] \setminus \{3\}$\\
\bottomrule
\end{tabular*}
\end{subtable}
\hfill
\begin{subtable}{0.22\textwidth}
\begin{tabular*}{\linewidth}[t]{@{\extracolsep{\fill}}ccccc|c}
\toprule
$x_1$ & $x_2$ & $x_3$ & $x_4$  & $x_5$ &  bears \\
\midrule
3 &3 &4 &5 &6  & 1, 2, 5\\
3 &3 &3 &4 &6  & 1, 2, 3 \\
3 &3 &3 &3 &5  & 2, 3, 4 \\
\hline
2 &3 &3 &3 &4  & 2, 3, 4 \\
1 &3 &3 &3 &3  & 3, 4, 5 \\
0 &2 &3 &3 &3  & 1, 4, 5 \\
\hline
0 &1 &2 &3 &3  & 1, 4, 5 \\
0 &0 &1 &3 &3  & 1, 2, 5\\
0 &0 &0 &2 &3  & 1, 2, 3 \\ 
\hline
0 &0 &0 &1 &2  & 1, 2, 3\\  
0 &0 &0 &0 &1 & 2, 3, 4 \\
-1 &0 &0 &0 &0 & 3, 4, 5 \\
\hline
-2 &-1 &0 &0 &0 & 3, 4, 5 \\
-3 &- 2& 0& 0& 0 & 1, 4, 5\\
-3 &-3 &-1 &0 &0 & 1, 2, 5\\
\hline
-3 &-3 &-2 &-1 &0 &1, 2, 5\\
\bottomrule
\end{tabular*}
\end{subtable}%
\hfill
\begin{subtable}{0.21\textwidth}
\begin{tabular*}{\linewidth}[t]{@{\extracolsep{\fill}}ccccc|c}
\toprule
$x_1$ & $x_2$ & $x_3$ & $x_4$ &  $x_5$ &  bears \\
\midrule
3 &3 &4 &5 &6 & 1, 2\\
3 &3 &3 &4 &5 & 2, 3\\
2 &3 &3 &3 &4 & 3, 4\\
\hline
1 &2 &3 &3 &3 & 4, 5 \\
0 &1 &2 &3 &3 & 1, 5\\
0 &0 &1 &2 &3 & 1, 2\\
\bottomrule
\end{tabular*}
\end{subtable}
\hfill
\begin{subtable}{0.2\textwidth}
\begin{tabular*}{\linewidth}[t]{@{\extracolsep{\fill}}ccccc|c} 
\toprule
$x_1$ & $x_2$ & $x_3$ & $x_4$  & $x_5$ &  bears \\
\midrule
3&  3&  4&  5&  6&  2\\
2&  3&  3&  4&  5&  3\\ 
1&  2&  3&  3&  4&  4\\ 
\hline
0&  1&  2&  3&  3&  1\\ 
0&  0&  1&  2&  2&  2\\ 
-1&  0&  0&  1&  1& 3\\ 
\hline
-2&  -1&  0&  0&  0&  5\\ 
-3&  -2&  -1&  -1&  0&  1\\ 
-3&  -3&  -2&  -2&  -1& 2\\
\hline
-4&  -3&  -3&  -3&  -2&  4\\
-5&  -4&  -4&  -3&  -3&  5 \\ 
-6&  -5&  -5&  -4&  -3&  1 \\
\hline
-6&  -6&  -6&  -5&  -4& 3\\ 
-7&  -7&  -6&  -6&  -5&  4\\
-8&  -8&  -7&  -6&  -6&  5\\
\hline
-9&  -9&  -8&  -7&  -6& 2\\
\bottomrule
\end{tabular*}
\end{subtable}
\hfill
\end{table}

\section{Proof of the main lemma} 
\label{s-proof-l-main}
Let us fix  $n, k, \ell,$ and $x = (x_1, \ldots, x_n)$. 
Recall that we always assume that the entries are monotone non-decreasing, $x_1 \leq \ldots \leq x_n$, and also that the inequality $n-k \leq m(x)$ holds, in other words,  all  $n-k$  bears are multiple of $\ell$. 

First, we consider in detail two special cases: $k=1$  and  $k = n-1$. 
Then we consider general case $1 < k < n$  skipping obvious details.  

\subsection{Case $k=1$}
In this case there are $n-k = n-1$ bears 
and only one bull, say, $i \in [n]$. 
Since $m(x) \geq n-k$, all entries of $x$, 
except maybe $x_i$, are multiples of $\ell$. 

We have to prove that in $\ell$ GM-moves 
the unique bull will be replaced by $i+1 \mod n$,  
in particular, by $n+1=1$ if $i=n$.
First, assume that $i < n$. 

Then, $x_{i+1} = \ldots = x_n = a \ell$  for some integer  $a$, 
while $x_{i-1} = (a-1) \ell$, since $range(x) \leq \ell$ 
and $i-1$ is a bear. 
Furthermore,  $(a-1) \ell < x_i \leq a \ell$, 
where the strict inequality holds by the tie-breaking rule.

Consider $\ell$ GM-moves from $x$ watching $x_{i+1}$. 
Being a unique bullish entry, 
$x_{i+1}$ is reduced by 1 with every GM-move 
until it reaches $(a-1) \ell$.  
This will certainly happen during the considered $\ell$  
GM-moves and, after this, by the tie-breaking rule, 
$i$  becomes a bear, $i+1$  becomes a unique bull, 
and both keep this status until 
the considered  $\ell$ GM-moves last. 
Thus, at the end  $i+1$  will be a unique bull. 

\medskip 

Now assume that a unique bull is $n$. 
There are two possible subcases: 
\begin{itemize}
\item[(i)] $x_n$  is a unique entry of $x$ which is
not a multiple of $\ell$; 
\item[(ii)] all entries of $x$  are multiples of $\ell$ 
but $x_n$ is a unique maximum.  
\end{itemize}
In fact, $x_n$  is the strict maximum in case (i) too. 
We have $x_n  = (a+1) \ell$, 
while $x_i = a \ell$  for all $i < n$.   
Consider $\ell$  GM-moves from $x$  watching  $x_1$. 
In both cases, (i) and (ii), 
$x_n$  is reduced by 1 by each move 
until  $x_n = x_{n-1}$  begins to hold. 

In case (ii) this happens in exactly $\ell$ GM-moves, 
all $n$  entries become equal, 
and $1$ becomes a unique bull, by the tie-breaking rule. 

In case (i), all $n$ entries become equal in less than $\ell$ GM-moves. 
Still $1$ becomes a unique bull and it keeps this status 
till the end of the considered $\ell$ moves, 
since the number of remaining GM-moves 
is strictly less than $\ell$ and, hence, 
$x_1$ does not have enough time to reach $(a-1) \ell$, 
becoming a bear.  

Thus, in both cases $1$ becomes the (unique) bull after $\ell$ 
GM-moves from $x$.

\subsection{Case $k=n-1$}
In this case there exists a unique bear 
$i \in [n]$, since $n-k=1$. 

We have to prove that in $\ell$ GM-moves 
this bear will be replaced by $i-1 \mod n$,  
in particular, by $1-1=0=n$ if $i=1$.
First, assume that $i > 1$. 

\medskip 

By the GM-rule, $x_i = a \ell$ for some integer $a$.  
If  $x_{i'} = a \ell$ then  $i' \leq i$, 
otherwise  $i' > i$  and, 
by the tie-breaking rule, $i'$ would be the bear, instead of $i$. 
The same happens if  $x_{i'}$  equals $(a-1) \ell$ 
for some $i' \neq i$.
Thus, $(a-1) \ell < x_{i'} \leq a \ell$  for all $i' \leq i$.

If $i' > i$ then $a \ell < x_{i'} \leq (a+1) \ell$. 
Indeed, $a \ell = x_i \leq x_{i'}$  results from monotonicity.

If equality holds, then by the tie-breaking rule,  
$i'$ is a bear rather than $i$.

Furthermore, $x_{i'} \leq (a+1) \ell$, 
since $range(x) \leq \ell$. 
For the same reasons, if  $x_{i'} = (a+1) \ell$ then  
$x_{i''} = a \ell$  for all $i'' \leq i$.   

Watch  $x_{i-1}$  for the next $\ell$  GM-moves. 
It is reduced by 1 with each move 
until it reaches $(a-1) \ell$, 
since all this time the bearish entry equals $a \ell$. 

Yet, in at most  $\ell$  GM-moves  
$x_{i-1}$  will be reduced to $(a-1) \ell$, 
since $x_{i-1} \leq x_i = a \ell$.  

Then $i-1$ becomes the bear and 
remains in this status 
until the considered  $\ell$  GM-moves end. 
This results from the following observations:    

If  $i' > i-1$  then $x_{i'}$  
cannot be reduced to $(a-1) \ell$ in $\ell$ GM-moves. 
Indeed, if  $i' > i$  then  $x_{i'} > a \ell$, as we know. 
If  $i' = i$  then  $x_{i'} = a \ell$, 
but still it cannot be reduced to $(a-1) \ell$
in $\ell$  GM-moves, 
because the first of them ``is lost", 
since  $i' = i$ is the bear in the beginning. 

If  $i' < i - 1$  then $x_{i'} \geq (a-1) \ell$, 
since $x_i = a \ell$ and  $range(x) \leq \ell$. 
If  $x_{i'} = (a-1) \ell$ 
then $i'$  is the bear rather than $i$. 
Hence,  $x_{i'} > (a-1) \ell$  and 
it cannot be reduced to 
$(a-2) \ell$  in $\ell$ GM-moves. 
Thus, after $\ell$  GM-moves $i-1$  will be the bear. 

\medskip 

Now assume that  $1$  is the bear. 
By the GM-rule,  $x_1 = a \ell$  for some integer~$a$.
Let us show that $a \ell < x_i \leq (a+1) \ell$ 
for all $i > 1$.
Indeed, inequality $x_i \geq a \ell$ holds, by monotonicity. 
If equality holds then, by the tie-braking rule, 
$i$  is the bear rather than 1. 
Hence, $x_i > a \ell$. 
Furthermore, $x_i \leq (a+1) \ell$, 
since $range(x) = x_n - x_1  \leq \ell$. 

Watch $x_n$ for the next $\ell$  GM-moves.
It is reduced by 1 with each move  
until $x_n$  reaches $a \ell$, 
since  all this time 
the bearish entry equals $a \ell$.  
Yet, in at most  $\ell$  GM-moves  
$x_n$  will be reduced to $a \ell$, 
since $x_n \leq (a+1) \ell$.  
Then $n$ becomes the bear and 
it remains in this status 
until the considered  $\ell$  GM-moves end 
since no entry can be reduced to $(a-1) \ell$.
Indeed, if  $i > 1$  then  $x_i > a \ell$, as we know, 
and $range(x) \leq \ell$.  
Although  $x_1 = a \ell$, 
but still it cannot be reduced to $(a-1) \ell$
in $\ell$  GM-moves, 
because the first of them ``is lost" 
since  $1$ was the bear in the beginning. 

\subsection{General case: $0 < k < n$}
Consider  $x = x^j$  for some $j \geq N(n,k,\ell,x^0)$. 
Then, $m(x) \geq n-k$, 
and hence, all $n-k$ bearish values of $x$  are multiples of $\ell$; 
also, $range(x) = x_n - x_1 \leq \ell$, 
and hence, among these  $n-k$ bearish values 
there are only 1 or 2 distinct. 

\medskip 

{\bf Subcase 1}  holds if and only if  
there exist integer $\mu$ and $\nu$  such that  
$$0 \leq \mu < \nu \leq n, \; \nu-\mu \geq n-k,$$  
$$(a-1) \ell < x_1 \leq \ldots \leq x_m < 
x_{\mu+1} = \ldots = x_\nu = a \ell 
< x_{\nu+1} \ldots \leq x_n$$ 
for some integer $a$.   
Then,  $n-k$ bulls and $k$ bears 
form the (cyclic) intervals   
\newline
$\{\nu-(n-k)+1, \ldots \nu\}$  and 
$\{1,\ldots,\nu_{n-k}\} \cup \{\nu+1,\ldots,n\}$, 
respectively, that is, bulls form an 
interval of type 1, 
while bears -- an interval of type 2.

Furthermore, 
$m(x) = \nu-\mu \geq n-k$.
Note that $x_i > (a-1) \ell$ 
for $i = 1, \ldots, \mu$, if any, 
and  $x_i < (a+1) \ell$ 
for $i = \nu+1, \ldots, n$, if any,  
since  $range(x) = x_n - x_1 \leq \ell$. 

\medskip 

In  $\ell$  GM-moves from $x$,  
the entries $x_1, \ldots, x_{\nu - (n-k)}$, if any,   
will be reduced to $(a-1) \ell$, 
then, by the tie-braking rule, 
at most  $n-k$  rightmost of them become bears.
Moreover, they remain in this status until 
the considered $\ell$  moves finished. 
Indeed, there is not enough time 
for  $x_i$  to reach  $(a-1) \ell$, since $i > \nu - (n-k)$. 
Although $x_{\nu - (n-k) + 1} = \ldots = x_\nu = a \ell$, 
but all these entries were bearish in the beginning, 
and hence, were not reduced by the first GM-move. 
Furthermore, $x_i > a \ell$  for $i > \nu$, 
and hence, these entries, if any, have not 
enough time to reach $(a-1) \ell$ either. 
However, they will be certainly reduced to $a \ell$, 
because $range(x) \leq \ell$.
Then, some number of the 
rightmost of these entries may become bearish,  
to make the total number of bears $n-k$, if necessary.
Furthermore, they remain in this status until 
the considered  $\ell$  moves finished. 

Let us summarize:  
After  $\ell$  GM-moves, 
the obtained $n-k$ bears form 
a (maybe, cyclic) interval that ends 
in $\nu - (n-k)$, 
while before these $\ell$ GM-moves,   
the  $n-k$  bears formed an interval 
beginning in $\nu-(n-k)+1$. 
Thus, considered $\ell$ GM-moves result 
in shifting the interval of $n-k$ bears 
by $n-k$ to the left, 
or equivalently, by  $k$ to the right.
The same is true for the complementary  
interval of $k$ bulls. 

\bigskip 

{\bf Subcase 2} holds if and only if  
there exist integer $\mu$ and $\nu$  such that  
$$0 < \mu < \nu < n, \;\; \nu - \mu \leq k,$$  
$$a \ell = x_1 = \ldots = x_\mu < 
x_{\mu+1} \leq \ldots \leq x_\nu < 
x_{\nu+1} = \ldots = x_n = (a+1) \ell$$ 
for some integer $a$.  
Bulls form the interval $\{\mu+1, \ldots, \mu+k\}$  
of size  $k$ and of type 1, 
while bears form an interval of type 1 or 2. 

In $\ell$ GM-moves the corresponding entries 
will be reduced to $a \ell$. 
In contrast, entries $x_{\mu+k+1}, \ldots, \mu_n$ 
have not enough time to be reduced to $a \ell$.
Although 
\newline 
$x_{\mu+k+1} = \ldots = \mu_n = (a+1) \ell$. 
yet, all these entries were bearish in the beginning, 
and hence, were not reduced by the first GM-move. 
Hence, after $\ell$ GM-moves, 
$\mu + k$  will become a bear, 
moreover, it will be an end of the 
interval formed by the  $k$ new bulls, 
while before these $\ell$ GM-moves, 
$\mu+k+1$  was the beginning of the 
(cyclic) interval of $k$ bears. 
Thus, considered $\ell$ GM-moves result 
in shifting the interval of $k$ bulls 
by $k$ to the left, 
or equivalently, by  $n-k$ to the right.
The same is true for the complementary 
interval of $n-k$  bears. 
\qed

\section{Proof of Proposition \ref{p1}}
\subsection{Case $k=n-1$; pursuit of the leaders}
\label{ss-pursuit} 
Recall that the  GM-sequence 
$S = S(n,k,\ell,x^0)$  
is partitioned 
into  $S_0 = \{0,1, \ldots, N\}$  and  
$S_\infty = \{N, N+1, \ldots\}$  such that  
$range(x^j) > \ell$  if $j < N$  and 
$range(x^j) \leq \ell$  if $j \geq N$.  
Here we prove that $N = N(n,k,\ell,x^0)$ 
grows as a polynomial and compute $N$ and $x^N$ 
in time linear in  $n,k,\ell$  and 
$\log(1 + range(x_0)) = \log(1 + x^0_n - x^0_1)$. 
   
Given $x = (x_1, \ldots x_n)$,  
an element $i \in [n]$ is called a {\em leader} 
if  $x_i - x_1 \leq \ell$. 
Let $\mu = \mu(x, \ell)$  be the maximal leader in $x$. 
Then, $[\mu] = \{1, \ldots, \mu\}$ is the set of leaders. 
Our algorithm is based on an iterative expansion of this set.

For a GM-sequence $S(n,k,\ell,x^0)$,  
obviously, $[\mu(n,k,\ell,x^0)] = [n]$  
if and only if  $range(x^0) \leq \ell$. 
This case is considered by the main lemma.  
Here we assume that $range(x^0) > \ell$, 
or in other words, the set of leaders in $x^0$ 
is a proper subset of $[n]$.  
Furthermore, it may happen that  $n-k > m(x_i^0)$. 
Yet, by Observation \ref{o3}, the inequality  
$n-k \leq m(x_i^j)$  will be achieved in at most 
$j = \frac{1}{k}[\sum_{i=1}^{m(x^0)-(n-k)}(x_i^0 \mod \ell)]-m(x^0)$  
GM-moves. 

Wlog, we can assume that it  
holds already for $j=0$, just to simplify the notation.  

Also wlog, we can assume that 
$[\mu(x)]$  will keep containing at least  $n-k$  bears  
until  $\mu+1$
(perhaps, together with $\mu+2$ and some other entries) 
catches up with the leaders. 
This will certainly happen and we can efficiently compute how soon. 

Consider the first $\mu \ell$ moves 
of the GM-sequence  $S$. 
By Theorem \ref{t-main},   
the entry $x^0_i$ will be reduced 
by  $(m -1) \ell$ if $i \leq m$,  
and by  $m \ell$ if  $i > m$. 
Thus, outsiders catch up with the leaders 
by  $\ell$  in every $m \ell$ GM-moves. 
Hence, in $m \ell \lfloor \frac{x_{m+1} - x_m}{\ell} \rfloor$  
GM-moves, the distance between them will be at most  $\ell$ 
and  in the next $\ell$ GM-moves, entry 
$m+1$, perhaps with some others, will join the leaders. 
Thus, we obtain a new initial vector and proceed iterating. 
In at most $n - m$  such iterations 
we obtain the last vector $x^N$  
such that $range(x^N) \leq \ell$, 
that is, the last vector of $S_0$. 

\subsection{General case: $0 < k < n$}
\label{ss-rs} 
Partition $[n]$  
by the first $n-k$ and last $k$ indices, 
$\{1, \dots, n\} = 
\{1, \dots, n-k\} \cup \{n-k+1, \dots, n\}$,  
that is, 
$[n] = [n-k] \cup ([n] \setminus [n-k])$.

By the GM-rule, 
$i \in [n-k]$   is a bull 
if it is not a multiple of $\ell$. 
Yet, in this case it becomes a
multiple of $\ell$ 
(and, hence, a bear) 
in at most $\ell-1$  GM-moves. 

If all $i \in [n-k]$ are multiples of $\ell$  
then  $i$ is a bear 
unless the next condition holds. 

\medskip 

(C)    
There exist integers $r, s$, and $a$  such that 
$0 < s \leq k, \; 0 < r \leq n-k$, and 
$x_i = a \ell$  if and only if 
$i \in I = [n-k-r+1, \dots, n-k+s]$, 
that is, $0 < n-k-r < i \leq n-k+s$. 

Then, by monotonicity, 
$x_i < a \ell$  if  $i \leq n-k-r$  and 
$x_i > a \ell$  if  $i > n-k+s$. 

\medskip 

Furthermore, if (C) holds then, by the tie-breaking rule, we have: 

\medskip 

(D) 
$i$  is a bull if  
$n-k-r < i \leq  n-k-r+s$  or $i > n-k+s$  and 
$i$  is a bear if  
$i \leq  n-k-r$  or  $n-k+s-r < i \leq  n-k+s$.

\medskip 

By counting,  
there are  $s + (k-s) = k$  bulls and 
$(n-k-r) + r = n-k$  bears. 

\medskip 

Note that the interval $I \subset [n]$ 
(from condition (C))  
contains $r+s$ successive entries 
and it intersects both intervals 
$[n-k]$  and  $[n] \setminus [n-k]$  
in  $r$  and  $s$  entries, respectively. 
If such $I$ does not exist and 
all $i \in [n-k]$  are multiples of $\ell$ 
then the pattern is simple:  
entries from $[n-k]$  are bears and 
from $[n] \setminus [n-k]$  are bulls.

Accordance to the GM-rule,  
the bearish entries stay, while 
all bullish are reduced by 1 
with each  GM-move.
Hence, the last $k$  entries 
are reduced until an interval  $I$  
satisfying  (C)  appears. 
The corresponding $j$  and $x^j$  
are easy to compute in time linear in 
$n,k,\ell, \log(1 + range(x^0))$, and $\log (1+j)$  
(rather than in $(x^0_n - x^0_1)$ and $j$).  

Then, $i$  becomes a bull if and only if   
$n-k-r < i \leq  n-k-r+s$  or $i > n-k+s$.  
The first $s$  bullish values    
will be reduced 
from $a \ell$  to $(a-1) \ell$ with $\ell$ GM-moves, 
thus turning these $s$ bulls into bears, then, the next  $s$, etc. 
In $\lceil r/s \rceil$  GM-moves 
$[n-k]$  becomes the set of all bears again. The rest of the interval $I$ will become bulls. 
Their bullish values will be reduced 
from $a \ell$ to $(a-1) \ell$ in the next $\ell$  GM-moves. 

These $\ell (1 + \lceil r/s \rceil)$ 
GM moves will be repeated several times. 
Each time all entries of the interval  $I$ 
are reduced by  $\ell$,  
taking values  $a \ell, (a-1) \ell, \ldots, (a-t) \ell$. 

\smallskip 

(i) Meanwhile, each $i'$ such that  
$i' < \min (i \mid i \in I)$  remains a bear, 
furthermore, $x_{i'}$ remains a constant multiple of $\ell$ 
and at most  $x_i \leq (a-t) \ell$.

\smallskip 

(ii) In contrast, each  $i''$  
such that  $i'' > \max (i \mid i \in I)$ remains a bull,  
its value $x_{i''}$  is reduced by 1 
with every GM-move but inequality 
$x_{i''} \leq (a-t) \ell$  holds.

\smallskip 

We proceed until the equality is achieved in (i) or (ii).
In both cases the interval $I$ will be expanded,  
while the number of distinct values 
taken by $x_i, i \in [n],$  reduced. 
Since this number is at most $n$, 
the procedure will be finished 
in at most $n$ steps and in $N$  GM-moves, 
resulting in a vector  $x^N$  
such that  $range(x^N) \leq \ell$. 
The algorithm sketched above implies 
that $N=N(n,k,\ell,x^0)$ is a polynomial 
and can be computed, together win $x^N$,
in time linear in 
$n,k,\ell, \log(1 + range(x^0))$, and $\log (1+j)$  
(rather than in $(x^0_n - x^0_1)$ and $j$).

\section{Finish line}
Applications to impartial games 
require the following generalization. 
Fix an integer  $d$  such that $0 < d \leq n$,   
introduce a finish level, say $f$, and assume that 
GM-sequence stops at $x^j$  as soon as 
at least  $d$  of its entries 
become of value at most  $f$.

For impartial games we will need only 
$d=2$ and $f=0$  but this case is not simpler 
than the general one. 

Suppose that $x^0$ already contains 
$d^0$  entries of value at most  $f$. 
If $d_0 \geq d$  then the GM-sequence 
is finished before it starts. 
Wlog, assume that $0 \leq d_0 < d$. 

\subsection{Case: $n = k+1$}
Suppose also that $x^0$ has $m_0$ leaders 
and proceed with the iterative process described 
in Section \ref{ss-pursuit}. 
We have to decide what will happen first: 

(i) $d$  entries become non-negative 
and the GM-sequence stops, or 

(ii) some new entries will join the leaders and $m$ increases. 

To do so, we compare 
$N_1 = m \ell \lfloor \frac{x_d}{\ell} \rfloor$ and 
$N_2 = m \ell \lfloor \frac{x_{m+1} - x_m}{\ell} \rfloor$. 

(j) If $N_1 \geq N_2$, we compute $x^{N_2}$ 
and proceed further with at most $\ell$  steps to finish 
the GM-sequence; 

(jj)  If $N_1 < N_2$, we expand the set of leaders and 
proceed with obtained  $m' > m$. 

Note that option (ii) disappears in case $range(x^0) \leq \ell$. 

\subsection{General case: $0 < k < n$}
Suppose that inequality 
$0 < x_i^0 < \ell$  holds 
for $m_0 = m_0(x^0)$ entries of $x^0$.  
If $n-k-m_0 \geq d-d_0$  
then at least $d$ non-positive entries 
appears in $x$  after at most $\ell-1$  GM-moves, 
thus terminating the procedure. 

Otherwise, if  $n-k-m_0 < d-d_0$,  
make $j_1 = \ell-1$  GM-moves and proceed further 
replacing  $x^0$ by $x^{j_1}$  and 
$d_0 = d_0(x^0)$  by $d_1 = d_0(x^{j_1})$.    
It is easily seen that 

\smallskip 

\noindent 
(E)  
$m(x^{j_1}) \geq n-k$, in other words, 
the first $n-k$  entries of $(x^{j_1})$ 
are multiples of $\ell$. 

\smallskip 

Phase I (which might be missing) ends 
as soon as property (E) holds, 
and the next Phase II lasts until it holds. 
Thus, during Phase II, 
$[n-k]$  are bears and 
$[n] \setminus [n-k]$ are bulls. 
The procedure can stop during phase II 
only if  $x_{n-k}^{j_1} \leq f < x_{n-k+1}^{j_1}$.
It happens as soon as $d - d_1$  bullish values 
become at most $f$. 
Clearly, the corresponding  $j$  and $x^j$ 
can be efficiently computed. 

However, phase II ends as soon 
as the interval $I \subset [n]$ 
intersecting both 
$[n-k]$  (in $r > 0$ entries) and  
$[n] \setminus [n-k]$ (in $s > 0$ entries) appears. 
Then, $s$ bears turn into bulls 
after the next GM-moves; see Section \ref{ss-rs}.

It may happen that only $d_2$ bullish values 
are reduced to $f$ or below 
during Phase II and  $d_1 + d_2 < d$. 
For example, $d_2 = 0$  (and $d_1 < d$)  
whenever $x_{n-k}^{j_1} > f$. 

In this case Phase III is required, 
which proceeds in accordance with Section \ref{ss-rs}:  
``$i$  becomes a bull if and only if   
$n-k-r < i \leq  n-k-r+s$  or $i > n-k+s$.  
\newline 
The first $s$  bullish values will be reduced 
from $a \ell$  to $(a-1) \ell$ with $\ell$ GM-moves, 
thus turning these $s$ bulls into bears, then, the next  $s$, etc. 
In $\lceil r/s \rceil$  GM-moves 
$[n-k]$  becomes the set of all bears again, while 
the rest of the interval $I$ become bulls. 
Their bullish values will be reduced 
from $a \ell$ to $(a-1) \ell$ in the next $\ell$ GM-moves.

These $\ell (1 + \lceil r/s \rceil)$ 
GM moves will be repeated several times. 
By each set, all entries of the interval  $I$ 
are reduced by  $\ell$,  
taking values  $a \ell, (a-1) \ell, \dots (a-t) \ell$."  

Then, follow (i) and (ii) of Section  \ref{ss-rs}, etc. 

Interval $I = I_0$  will be successively replaced by 
$I_1, \dots I_\tau$, where  $\tau \leq n$. 
\newline 
Respectively, Phase III is partitioned into 
sub-phases III-1 $\dots$ III-$\tau$.  
For each sub-phase III-$t$, $0 \leq t \leq \tau$,  
one can efficiently compute 
the number $d_t$ of entries that become at most $f$. 
The whole procedure stops as soon as 
$d_0 + d_1 + \dots + d_t \geq d$. 

Yet, Phase III stops as soon as $range(x^j) \leq \ell$ holds, 
which may happen earlier. 
As we know, the latter inequality is respected by the GM-moves. 

Then, Phase IV begins, 
during which $x$ is changed ``periodically", 
in accordance with formula (\ref{eq-main}), 
Lemma \ref{l-main},  and Theorem \ref{t-main}.

We have to determine $j$ and $x^j$ such that 
$d_0 + d_1 + \dots + d_j \geq d$,
where $d_{\nu}$ is the number of entries 
$x_i^\nu \leq f$ for $i \in [n]$ and $0 \leq \nu \leq j$.

Both   $j$ and $x^j$  can efficiently computed 
in each Phase III-t or IV. 
We leave this simple but tedious case analysis 
to the careful reader.

\section{Remoteness function of exact slow NIM; 
case $k=n-1,\ell=2,f=0$}
\label{s-games}
We assume that the reader is familiar 
with basic concepts of impartial game theory
(see e.g., \cite{ANW07,BCG01-04} for an introduction) 
and also with the recent paper \cite{GMMV23}, 
where the game NIM$(n,=k)$ of exact slow NIM,   
was analyzed for the case $n = k+1$.  
Its mis\`ere version was considered in \cite{GMMN23}.
Games NIM$(4,=2)$ and NIM$(5,=2)$ 
were recently considered in \cite{GMN23}.  
It appears that NIM$(n,=k)$ is closely related 
the screw discrete dynamical system 
with parameters $n$ and $k$. 
Furthermore, $f=0$, $\ell=2$, and $d=n-k+1$. 
In particular, $d =2$  for $n = k+1$.

\subsection{Exact Slow NIM} 
\label{ss00}
Game Exact Slow NIM  
was introduced in \cite{GH15} as follows: 
Given two integers  $n$ and $k$  
such that  $0 < k \leq n$  and  
$n$  piles containing  $x_1, \ldots, x_n$  stones.  
By one move it is allowed to reduce each of any  $k$  piles 
by exactly one stone. 
Two players alternate turns. 
The player who has to move but cannot is the loser 
in the normal version of the game and 
the winner of its mis\`ere version. 
Obviously, $x$ is a terminal position 
if and only if less than $k$ entries of $x$ are positive. 
In \cite{GH15}, this game was denoted  NIM$^1_=(n,k)$.
Here we will simplify this notation to  NIM$(n,=k)$.

It is easily seen that  NIM$(n,=k)$  
is partitioned into  $k$  disjoint subgames, 
since $x_1 + \ldots + x_n \mod k$  is an invariant, 
it is not changed with moves. 

Note also that $x = (x_1, \ldots, x_n)$ 
is a P-position whenever  $x_i$  is even for all $i \in [n]$. 
Indeed, the second player can keep this property 
just always repeating the move of the first player. 

Game NIM$(n,=k)$  is trivial if  $k = 1$  or  $k = n$. 
In the first case it ends after $x_1 + \ldots + x_n$ moves  
and in the second one---after  $\min(x_1, \ldots, x_n)$ 
moves. In both cases nothing depends on players' skills. 
All other cases are more complicated. 

The game was solved for  $k=2$  and  $n \leq 6$.
In \cite{GHHC20}, an explicit formula 
for the Sprague-Grundy 
(SG) function was found for  $n \leq 4$ and $k = 2$, 
for both the normal and mis\`ere versions. 
This formula allows us to compute the SG function in linear time. 
Then, in \cite{CGKPV21} the P-positions 
of the normal version were found for  $n \leq 6$ and $k = 2$.  
For the subgame with even  $x_1 + \ldots + x_n$ 
a simple formula for the P-positions was obtained, 
which allows to verify in linear time 
if  $x$ is a P-position and, if not, 
to find a move from it to a P-position.  
The subgame with odd  $x_1 + \ldots + x_n$ is more difficult. 
Still a (more sophisticated)  formula 
for the P-positions was found,  
providing a linear time recognition algorithm.

\subsection{Case $n = k+1$}
\label{ss01}
In \cite{GMMV23} the normal version of the game  
was solved in case $n = k+1$ 
by the following simple rule: 
\begin{itemize}
\item[(o)] if all piles are odd, 
keep a largest one and reduce all other;
\item[(e)] if there exist even piles, 
keep a smallest one of them and reduce all other.
\end{itemize}

This rule was called the {\em M-rule} in \cite{GMMV23}. 
Obviously, it coincides with the GM-rule,  
if we restrict ourselves by parameters  
$n=k+1$  and  $\ell = 2$; 
furthermore, the terminating parameters are
$f=0$  and $d=n-k+1=2$.
Obviously, $x$ remains non-negative 
during the play, which terminates 
as soon as at least two entries of $x$ 
become 0. 

Denote by $\M(x)$  the number of moves 
from  $x$  to a terminal position, 
assuming that both players follow the M-rule. 
In \cite{GMMV23} it was proven that  $\M = \R$, where 
$\R$  is the {\em remoteness function} 
introduced by Smith in \cite{Smi66}. 
Thus, M-rule solves the game and, moreover, 
it allows to win as fast as possible in an N-position and 
to resist as long as possible in a P-position.

A polynomial algorithm computing  $\M = \R$  
(and in particular, the P-positions) was given, 
even if $n$ is a part of the input and 
integers are presented in binary form.
Results of the present paper, restricted to $\ell = 2$, 
provide a simpler algorithms.

Let us also note that 
an explicit formula for the P-positions is known 
only for $n \leq 4$  
and already for $n=3$ it is pretty complicated 
\cite{GHHC20}\cite[Appendix]{GMMV23}.

\subsection{Games NIM$(4,2)$ and NIM$(5,2)$} 
It was shown in \cite{GMN23} 
that the M-rule is also applicable for computing 
the remoteness function of the mentioned two games. 
Equality  $\R = \M$  holds frequently.  
The exceptions are sparse, 
have regular patterns, and described by simple closed formulas.
However, these formulas are not proven. 
Recall that $f=0$ and $d=n-k+1$. 

\subsection{Related versions of NIM} 
By definition, the present game NIM$(n,k)$ 
is the exact slow version  
of the famous Moore's NIM$_k$ \cite{Moo910}.  
In the latter game a player, by one move, 
reduces arbitrarily 
(not necessarily by one stone) 
at most $k$ piles from $n$. 

The case  $k=1$  corresponds to the classical NIM  
whose P-position was found by Bouton 
\cite{Bou901}  for both the normal and mis\`ere versions. 

\begin{remark} 
Actually, the Sprague-Grundy (SG) values
of NIM were also computed in Bouton's paper, 
although were not defined explicitly in general. This 
was done later by Sprague \cite{Spr36} 
and Grundy \cite{Gru39} for 
arbitrary disjunctive compounds 
of impartial games; see also \cite{Con76,Smi66}. 

In fact, the concept of a P-position 
was also introduced by Bouton in \cite{Bou901}, 
but only for the (acyclic) digraph of NIM, 
not for all impartial games. 
In its turn, this is a special case 
of the concept of a kernel, 
which was introduced for arbitrary digraphs 
by von Neumann and Morgenstern \cite{NM44}. 

Also the mis\`ere version 
was introduced by Bouton in \cite{Bou901},  
but only for NIM, not for all impartial games. 
The latter was done by Grundy and Smith 
\cite{GS56}; see also \cite{Con76,Gur07,Gur07a,GH18,Smi66}.
\end{remark} 

Moore \cite{Moo910} obtained an elegant explicit formula 
for the P-positions of NIM$_k$ 
generalizing the Bouton's case  $k=1$. 
Even more generally, the positions 
of the SG-values 0 and 1  were efficiently characterized by 
Jenkins and Mayberry \cite{JM80}; see also 
\cite[Section 4]{BGHMM17}.  
Also in \cite{JM80};
the SG function of NIM$_k$  
was computed explicitly for the case  $n = k+1$ 
(in addition to the case $k=1$).  
In general, no explicit formula, 
nor even a polynomial algorithm, 
computing the SG-values 
(larger that 1)  is known. 
The smallest open case: 
2-values for  $n=4$ and $k=2$. 

\medskip 

The remoteness function of  NIM$_k$ was recently studied 
in \cite{BGMV23}. 

\medskip 

Let us also mention 
the exact (but not slow) game NIM$^1_=(n,k)$  
\cite{BGHM15,BGHMM17} 
in which exactly $k$ from $n$ piles are reduced 
(by an arbitrary number of stones) in a move. 
The SG-function was efficiently computed 
in \cite{BGHM15,BGHMM17} for $n \leq 2k$. 
Otherwise, even a polynomial algorithm 
looking for the P-positions is not known
(unless $k=1$, of course). 
The smallest open case is  $n=5$  and  $k=2$.

\subsection{Game NIM$(n,n-1)$  for $\ell$ players} 
The exact slow NIM$(n,k)$, with $k = n-1$  
can be played by $\ell$  players  
as follows. 
A position is a non-negative $n$-vector  $x$. 
Players make moves in a given cyclic order. 
By one move a player can choose an arbirary 
entry and keep it unchanged 
reducing $n-1$ remaining entries by 1, 
provided they are positive. 
If  $x$  has at least two non-positive entries, 
it is called a {\em terminal} position, 
the game is over, and the player 
who has to move in $x$ (but cannot) is a loser, 
while  $\ell-1$  other players are winners. 
The payoffs are defined below. 
A sequence of moves 
from the initial position  $x^0$ 
to a terminal one is called a {\em play}. 
Let  $L(P)$  denote the length, 
that is, the number of moves, of a play  $P$. 
Choose a large constant $C$; 
it should be larger than the length 
of any play from  $x^0$.  
Then the loser pays $C-L(P)$ 
and $\ell-1$  winners share this amount, 
gaining $\frac{C - L(P)}{\ell-1}$ each. 
The GM-rule defines the unique strategy for each player.
If a position  $x$  has no entries which are multiples of $\ell$, 
the GM-move in $x$ keeps the largest entry and reduces 
the smaller $n-1$ by 1 each. 

\begin{conjecture}
\label{con-1} 
The set of $\ell$ GM-strategies form a (uniform) Nash equilibrium. 
\end{conjecture}

By definition, the GM-strategies are uniform, 
that is, independent of  $x^0$. 
Hence, a  Nash equilibrium in GM-strategies  
is uniform, if exists. 

For $\ell = 2$  this conjecture immediately 
follows from the results of \cite{GMMV23}. 

\section{Concluding remarks 
and open problems} 
The GM-$(n,k,\ell)$-rule defines a deterministic dynamic system. 
It is well known that such systems may demonstrate 
a ``chaotic behavior". 
Yet, the system is very simple in the considered case, 
as shown by Lemma \ref{l-main} and Theorem \ref{t-main}.  
After the first $N$, every next $p$  GM-moves  
reduce all entries of $x$  by the same constant $\delta$; 
see formula (\ref{eq-main}).  
Furthermore, $N$ and $x^N$  can be efficiently computed.

Given a non-negative integer vector $x$, and  
integers  $\ell \geq 2, d \geq 1$, and $f$,  
we can determine, in time linear in 
$n,k,\ell,d$, $\log(1+x_n-x_1)$ and $\log(1+|f|)$,   
how many GM-moves are required to get 
$x'$  with at least $d$ entries of values at most $f$. 
In case $\ell=d=2$ and $f=0$, this number 
is the value of the remoteness function  $\M(x) = \R(x)$ 
of the slow NIM game NIM$(n,k)$  with $k=n-1$.   
This game can be extended to the case of $\ell \geq 2$ players; 
see Conjecture \ref{con-1}.

Interestingly the obtained optimal GM-strategies 
are uniformly optimal in the following sense. 
Let us replace  0  by an arbitrary 
integer $c$, positive or negative, and  
require that at least  $n - d$ entries of  $x$ are at least $c$. 
The game is over as soon as at least 
$d$  entries become at most  $c$. 
In the obtained game 
the GM-strategies are optimal and 
the same for any {\em even} $c$. 
Which other impartial games admit such 
uniformly optimal strategies? 

A more general open question:     
which other discrete dynamic systems 
are related to impartial games;  
in particular, it remains open 
already for the systems defined 
by the GM-rule with $\ell > 2$. 

 \bigskip 

 {\bf A $d$-dimensional generalization} 
 \newline
 Consider a $d$-dimensional table 
 $N = [n_1] \times \dots \times [n_d]$ 
 and a mapping  $x : N \to \mathbb{Z}$, 
 where $\mathbb{Z}$ is the set of integers. 
 We assume that  $x$  is entry-wise 
 monotone non-decreasing, that is, 
 $x(i) \leq x(i')$  for all 
 $i = (i_1, \dots, i_d)$, 
 $i' = (i'_1, \dots, i'_d)$, 
 and  $\alpha_0 \in [d]$  such that 
$i_\alpha, i'_\alpha \in [n_\alpha]$  and 
 $i_\alpha = i'_\alpha$  for all  $\alpha \in [d]$, 
 except only one, $\alpha_0 \in [d]$,  
 for which  $i_{\alpha_0} \leq i'_{\alpha_0}$. 
 Also consider the standard lexicographical order over $N$. 
 
Fix integer $k$ and $\ell$  such that 
$0 < k < |N| = n = n_1 \times \dots \times n_d$  and $1 < \ell$.
Denote by $m = m(x)$  the number of entries of 
$i \in N$  such that $x(i)$ is a multiple of $\ell$. 

A unique GM-move from $x$ is defined as follows. 
If $m(x) < n-k$  
then choose all vectors $i \in N$  such that 
$x(i)$ is multiple of $\ell$  and 
add  $k$  more entries of  $N$ arbitrarily, 
for example, the lexicographical largest ones. 
If $m(x) \geq n-k$ then choose 
$n-k$  lexicographical largest entries of $N$ 
multiple of $\ell$. 
(This is the tie-breaking rule.)  
By the GM-move from $i$, 
the chosen $n-k$ entries of $N$ (bears) keep their values, 
while the remaining $k$  (bulls) are reduced by 1.  

It is not difficult to verify that a GM-move 
from $x$ respects 
the entry-wise monotonicity of $x$, 
as well as the inequalities 
$m(x) \geq n-k$ and $range(x) \leq \ell$. 

Important open questions are:  
\newline
(i) Is the generalized GM-sequence $S$ quasi-periodical?  
\newline
(ii) If yes, what is its period $p$ and step $\delta$? 
\newline
(iii) Can vector $x^j$ from $S$  be computed in time 
polynomial in $n, k, \ell$, and  $\log(1+j)$?

 \bigskip 

\noindent 
{\bf Acknowledgements.}
This research was supported by Russian Science Foundation, 
grant  20-11-20203, https://rscf.ru/en/project/20-11-20203/.

\end{document}